\documentclass[12pt]{amsart}
\usepackage[all]{xy}
\usepackage{amssymb}
\calclayout
\newtheorem{dummy}{anything}[section] 
\newtheorem{theorem}[dummy]{Theorem}
\newtheorem*{thma}{Theorem A}

\newtheorem{lemma}[dummy]{Lemma} 
 
\newtheorem{corollary}[dummy]{Corollary}
 
\theoremstyle{definition}

 \newtheorem{example}[dummy]{Example}
 
 \newtheorem{remark}[dummy]{Remark}

\newcommand{\cA}{\mathcal A}
\newcommand{\cG}{\mathcal G}
\newcommand{\cH}{\mathcal H}
\newcommand{\cM}{\mathcal M}

\newcommand{\cD}{\mathcal D}

\newcommand{\cF}{\mathcal F}
\newcommand{\cK}{\mathcal K}
\newcommand{\cI}{\mathcal I}

\newcommand{\bZ}{\mathbf Z}
\newcommand{\bQ}{\mathbf Q}

\newcommand{\bC}{\mathbb C}

\newcommand{\cy}[1]{\bZ/{#1}}
\newcommand{\Zhat}{\widehat \bZ}

\newcommand{\mmatrix}[4]{\left (\vcenter
{\xymatrix@C-2pc@R-2pc{#1&#2\\#3&#4} }
\right )}

\DeclareMathOperator{\wh}{Wh}
\DeclareMathOperator{\Mod}{mod}

\DeclareMathOperator{\Image}{Im}
\DeclareMathOperator{\Ind}{Ind}
\DeclareMathOperator{\Res}{Res}
\newcommand{\newmod}[2]{#1\, (\Mod{#2})}

\begin{document}
\title[Free actions on Products of Spheres]
{Some Examples of Free actions on  Products\\ of Spheres}
\author{Ian Hambleton}
\address{Department of Mathematics \& Statistics
\newline\indent
McMaster University
\newline\indent
Hamilton, ON L8S 4K1, Canada}
\email{ian@math.mcmaster.ca}
\date{Jan. 25, 2006}
\thanks{\hskip -11pt  Research partially supported by NSERC Discovery Grant A4000. }

\begin{abstract} If $G_1$ and $G_2$ are finite groups with
periodic Tate cohomology, then $G_1\times G_2$
acts freely and smoothly  on some product $S^n \times S^n$.
\end{abstract}
\maketitle
\section*{Introduction}
The study of free finite group actions on products of spheres is
a natural continuation of the spherical space form problem
 \cite{madsen-thomas-wall1}. 
 In this paper, we show that certain products of finite groups do act freely and
\emph{smoothly} on a product of spheres, even though
the individual factors can't act freely and smoothly (or even topologically) on a single sphere. 
This verifies a conjecture of  Elliott Stein \cite{stein1}. The method
involves a detailed analysis of the product formulas in surgery theory, and a refinement of Dress induction for surgery obstructions.

\smallskip
If a finite group $G$ acts freely  on $S^n$, then
(i) every abelian subgroup of
$G$  is cyclic,
and (ii) every element of order $2$  is central. In
 \cite{madsen-thomas-wall1},
Madsen, Thomas and Wall proved that these conditions are
sufficient to imply the existence of a free topological action on some sphere.
Actually, these two conditions have a very different character.
By the work of P. A. Smith and R. Swan \cite{swan1}, 
condition (i) is necessary and sufficient for a free simplicial action of $G$ on a finite-dimensional \emph{simplicial complex}, which is homotopy equivalent to a sphere. 
The finite groups $G$ satisfying condition (i) are exactly the groups
with periodic Tate cohomology, or equivalently those for which
every subgroup of order $p^2$, $p$ prime, is cyclic
(the $p^2$-conditions).
On the other hand, Milnor \cite{milnor1} proved that condition (ii) is necessary for a free $G$-action by homeomorphisms
on any closed, topological \emph{manifold} which has the mod $2$
homology of a sphere. The groups with periodic cohomology satisfying condition (ii) are just those
which have no dihedral subgroups, or equivalently
those for which every subgroup of order $2p$, $p$ prime, is cyclic
(the $2p$-conditions). Milnor's result shows for example that the periodic dihedral groups do not act topologically
on $S^n$, although they do act simplicially on a finite
complex homotopy equivalent to $S^n$.
 
\smallskip
For free finite group actions on a  product of spheres,
 the analogue of
condition (i) was suggested  by P. Conner  \cite{conner1}:
 if $G$
acts freely on a $k$-fold product of spheres 
$(S^n)^k:=S^n \times \dots \times S^n$,
 is every abelian subgroup of $G$   generated by at most
$k$ elements ? Conner proved this statement for $k=2$,
and a lot of work \cite{oliver1}, \cite{carlsson1}, \cite{adem-browder},
\cite{adem-smith} has been done
to determine what additional conditions are necessary to produce
free simplicial actions on a finite-dimensional simplicial complex homotopy 
equivalent to a product of spheres. The picture is now almost
completely clarified, at least for elementary  abelian groups
and spheres of equal dimension:
Adem and Browder \cite{adem-browder} and Carlsson 
\cite{carlsson1} showed that  $G = (\bZ/p)^r$ acts freely on $(S^n)^k$, for
$p$ a prime, implies $r\leq k$ provided that $n \neq 1,3,7$ in the case
$p=2$ (the restriction $n\neq 1$ for $p=2$ was recently removed
in \cite{yalcin1}). The same result is conjectured to hold for
finite-dimensional $G$-$CW$ complexes homotopy equivalent
to a product of spheres of  possibly \emph{unequal} dimensions
(see \cite{adem-davis-unlu} for some recent progress).

Much less seems to be known at present about the additional conditions
needed to produce free actions by homeomorphisms 
or diffeomorphisms on the closed
manifolds
$(S^n)^k$ for $k >1$.
Let $D_q$ denote
the dihedral group of order $2q$, with $q$ an odd prime.
Elliott Stein \cite{stein1} proved that,  for every $n = 4j +3$ and
any $k\geq 2$, there exist free, orientation-preserving piece-wise
linear actions of $(D_q)^k$ on $(S^n)^k$.
Many of these  actions are smoothable.

These examples show that a direct generalization of Milnor's
condition (ii) is \emph{not} necessary for actions on products
of spheres.  In this paper we verify a conjecture of Stein's:
\begin{thma} 
 If $G_1$ and $G_2$ are finite groups with
periodic Tate cohomology, then $G_1\times G_2$
acts freely and smoothly  on some product $S^n \times S^n$.
\end{thma}
The techniques used to prove this statement also show that
any product 
of periodic groups $G_1 \times \dots \times G_k$, with
 $k>1$, 
acts freely and topologically  on  $(S^n)^k$ for some $n$. These are smooth actions if $k\neq 3$ (see Theorem \ref{thm: moreex}).
Of course there are groups $G$ satisfying Conner's condition
which are not the direct product of periodic groups, so these examples
are just the simplest case. The surgery techniques need to be
developed further to study more general groups. I would like to 
thank Alejandro Adem for reminding me about E. Stein's paper and this open question, and Andrew Ranicki for a correction to Lemma 
\ref{lem: nattrans}. I would also like to thank the referee for useful comments.

\section{Hyperelementary induction}
We need a refinement of Dress's fundamental work on induction, following
\cite[\S 1]{htaylor1} and \cite{htw4}. In 
 \cite{dress1} and
\cite{dress2}   we are given a Mackey
functor ${\cM}$ and a \emph{family} ${\cH}$ of subgroups of $G$,  which by definition is a collection of subgroups closed
under conjugation and taking further subgroups. 
An important example is the family of hyperelementary subgroups.

One can then form what Dress calls an Amitsur complex: 
this is a chain complex
$$\xymatrix{\cM (G)\ar[r]^(0.4){\partial_0}&\bigoplus_{H\in{\cH}}{\cM}(H)
\ar[r]^(0.7){\partial_1}&\dots}$$
where the higher terms are  explicitly described sums of 
${\cM}$ applied to elements of ${\cH}$.
The boundary map $\partial_0$ is the sum of restriction
maps and the higher $\partial_i$ are just sums and differences
of restriction maps.
There is a second Amitsur complex defined using induction maps
for which the boundary maps go the other direction.

Dress further assumes that some Green ring, say ${\cG}$,
acts on ${\cM}$.
Write 
$$\delta^{{\cH}}_{{\cG}}\colon \bigoplus_{H\in{\cH}}
{\cG}(H)\to {\cG}(G)$$
for the sum of the induction maps.
\begin{theorem}

If there exists $y\in \bigoplus_{H\in{\cH}}{\cG}(H)$
such that $\delta^{{\cH}}_{{\cG}}(y)=1\in{\cG}(G)$,
then both Amitsur complexes for ${\cM}$ are contractible.
\end{theorem}

\begin{remark}
The conclusion can be expressed as
$${\cM}(G)=
\lim\limits_{\stackrel{\longleftarrow}{\cH}}
{\cM}(H)\qquad \text{or}\qquad  {\cM}(G)=
\lim\limits_{\stackrel{\longrightarrow}{\cH}}
{\cM}(H)$$ where the first limit made up of restrictions
and the second of inductions. In these cases we say that
$\cM$ is $\cH$-detected or $\cH$-generated. If both
Amitsur complexes are contractible we say that $\cM$ is
$\cH$-computable.
The result above follows from \cite[Prop. 1.2, p.305]{dress2}
and the remark just above \cite[Prop. 1.3, p.190]{dress2}.
\end{remark}

In order to use this strong theorem, it is necesary to find a 
computable  Green ring which acts on a given Mackey functor.
The Burnside ring of finite $G$-sets 
is a Green ring which  acts on any
Mackey functor, via the formula $$[G/H]\cdot x = \Ind_H^G(\Res_G^H(x)), $$
 but it satisfies Dress's
condition on $\delta^{{\cH}}_{{\cG}}$
if and only if $G\in{\cH}$.
\begin{example}\label{hyperL}
The Dress ring $GU(G,\bZ)$ is
 $\cH$-computable for $\cH$ the hyperelementary
family (see \cite[Thm. 3]{dress2}). Since there is a homomorphism given
by tensor product from $GU(G,\bZ)$ to the automorphisms of
$(RG,\omega)$-Morita in $(R,-)$-Morita 
(see \cite[1.B.3]{htw3} for the definitions), it follows that
any functor out of  $(RG,\omega)$-Morita is also hyperelementary
computable. For this paper, the main examples will be
the quadratic and symmetric $L$-groups $L_n(\bZ G)$  and
$L^n(\bZ G)$, with various torsion decorations, as defined by Ranicki \cite{ra10}.
These are functors out of $(\bZ G, \omega)$-Morita, and hence are
hyperelementary computable (compare \cite[Thm. 1]{dress2} which
covers the quadratic $L$-groups at least).
Since the $2$-localization map $L_n(\bZ G)\to L_n(\bZ G)\otimes
\bZ_{(2)}$ is injective for $G$ finite, the $2$-hyperelementary
family suffices for deciding when a surgery obstruction is zero.
\end{example}

For each prime $p$, and for each subgroup $K$, we denote by $\cF_p(K)$ the set of subgroups
of $H \subset G$ where $K$ is a normal subgroup of $H$ and
$H/K$ is a $p$-group. A family $\cF$ is called $hyper_p$-\emph{closed} if $\cF_p(K) \subseteq \cF$ for all $K \in \cF$.

Let $X$ be a finite $G$-set. The minimal representative $G$-set for
the family of subgroups $\{\cF_p(K): K \in Iso(X)\}$ will 
be denoted $hyper_p$-$X$. This construction is due to Dress
\cite[p. 307]{dress2}. The trivial $G$-set is denoted $\bullet$.
One of Dress's main results is the following:
\begin{theorem}[{\cite[p. 207]{dress1}}] Let $\cG$ be a Green ring. 
For any prime $p$, and any finite $G$-set $Y$, let 
 $\cK(Y):=\ker\left ( \cG(\bullet) \otimes \bZ_{(p)} \to \cG(Y) \otimes \bZ_{(p)}\right )$ and let
$\cI(Y):= \Image\left ( \cG(hyper_p\text{-}Y)\otimes \bZ_{(p)} \to \cG(\bullet)\otimes \bZ_{(p)}
\right )$. Then $\cG(\bullet) \otimes \bZ_{(p)} = \cK(Y) + \cI(Y)$.
\end{theorem}
If $X$ is a finite $G$-set, 
we will use the notation $\langle X\rangle$ for the equivalence class
of $X$ in the Dress category $\cD(G)$.
One useful consequence is:
\begin{lemma}[\cite{htw4}] Let $\cG_0$ be a sub-Green ring of $\cG_1$. For any prime $p$, and any finite $G$-set
$X$ with $\langle X \rangle = \langle hyper_p\text{-}X\rangle$, then
the natural map $\cG_0(X) \otimes \bZ_{(p)} \to \cG_0(\bullet) \otimes \bZ_{(p)}$ is surjective if and only if 
$\cG_1(X) \otimes \bZ_{(p)} \to \cG_1(\bullet) \otimes \bZ_{(p)}$
is surjective.
\end{lemma}
\begin{proof}
For any Green ring $\cG$ and any finite $G$-set $Y$, 
the image of $\cG(Y) \otimes \bZ_{(p)}$ in $ \cG(\bullet)$ is an ideal.
Hence either map is onto if and only if $1_{\cG_i(\bullet)}$
is in the image. Since $1_{\cG_0(\bullet)}$ goes to
$1_{\cG_1(\bullet)}$, this proves the first implication.

For the converse, notice that the assumption 
$\cG_1(X) \otimes \bZ_{(p)} \to \cG_1(\bullet) \otimes \bZ_{(p)}$
is surjective implies that the Amitsur complex is contractible for the restriction maps
induced by the transformation $X \to \bullet$. In particular,
$\cG_1(\bullet) \otimes \bZ_{(p)} \to \cG_1(X) \otimes \bZ_{(p)}$
is injective. Therefore $\cG_0(\bullet) \otimes \bZ_{(p)} \to \cG_0(X) \otimes \bZ_{(p)}$ is injective, and from Dress's Theorem we conclude that
$\cG_0(hyper_p\text{-}X) \otimes \bZ_{(p)} \to \cG_0(\bullet) \otimes \bZ_{(p)}$ is surjective. 
\end{proof}
Suppose that $\cG$ is a Green ring which acts on a Mackey
functor $\cM$. 
For many applications of induction theory, 
the ``best" Green ring for $\cM$ is the \emph{Burnside quotient Green
ring} ${\cA}_{{\cG}}$,  defined as
the image of the Burnside ring in ${\cG}$.
This  is a Green
ring which acts on ${\cM}$, and by construction 
${\cA}_{{\cG}}$  is a sub-Green ring of $\cG$. In particular,
the natural map ${\cA}_{{\cG}} \to \cG$ is an injection.

We say that a finite $G$-set $X$ is a \emph{Dress generating
set} for a Green ring $\cG$, provided that, for each prime $p$, 
$\cG(hyper_p\text{-}X) \otimes \bZ_{(p)} \to \cG(\bullet) \otimes \bZ_{(p)}$ is surjective. 

\begin{theorem}
A finite $G$-set $X$ is a Dress generating set for a Green ring
$\cG$ if and only if it is a Dress generating set for the Burnside
quotient Green ring $\cA_{\cG}$.
\end{theorem}
We can translate this into a computability statement as follows:
\begin{corollary} Let $p$ be a prime and $\cG$ be a Green ring. Suppose that 
$\cF$ is a $hyper_p$-closed family of subgroups of $G$.  Then $\cG\otimes \bZ_{(p)}$ is $\cF$-computable
if and only if $\cA_{\cG}\otimes \bZ_{(p)}$ is $\cF$-computable.
\end{corollary}
The advantage of ${\cA}_{\cG}$ over ${\cG}$
is that ${\cA}_{\cG}$ acts on Mackey functors
which are subfunctors or quotient functors of ${\cM}$
but ${\cG}$ does not in general.
For example, ${\cG}$ never acts on ${\cA}_{\cG}$
unless they are equal. 
We next point out another good  feature
 of the Burnside quotient Green ring.

\begin{theorem}\label{thm: computable}
 Suppose that $\cG$ is a Green ring which
acts on a Mackey functor $\cM$,  and 
$\cF$ is a $hyper_p$-closed family of subgroups of $G$. If $\cG\otimes \bZ_{(p)}$ is $\cF$-computable, then  every
$x\in \cM(G)\otimes \bZ_{(p)}$ can be written as
$$x = \sum_{H \in \cF} a_H \Ind_H^G(\Res_G^H(x))$$
for some coefficients $a_H \in  \bZ_{(p)}$.
\end{theorem}
\begin{proof} 
Since $\cG\otimes \bZ_{(p)}$ is $\cF$-computable, we know that $\cA_{\cG}\otimes \bZ_{(p)}$ is
also $\cF$-computable. Therefore, we can write 
$1 = \sum_{K\in \cF} b_K \Ind_K^G(y_K)$, for some $y_K \in 
\cA_{\cG}(K)\otimes \bZ_{(p)}$ and $b_K\in\bZ_{(p)}$. For  any $x
\in \cM(G)\otimes \bZ_{(p)}$ we now have the formula
$$x = 1\cdot x = 
\sum_{K\in \cF} b_K \Ind_K^G(y_K\cdot \Res_G^K(x))$$
But each $y_K \in \cA_{\cG}(K)\otimes \bZ_{(p)}$ can be represented by 
a sum $\sum c_{KH}[K/H]$, with $c_{KH}\in  \bZ_{(p)}$, under the surjection $\cA(K) \to \cA_{\cG}(K)$. Therefore
$$
\begin{array}{ll}
x & = \sum_{K\in \cF} b_K
\sum_{H \subseteq K} c_{KH} \Ind_K^G([K/H]\cdot \Res_G^K(x))\\
&\\
& = \sum_{K\in \cF} b_K\sum_{H \subseteq K} c_{KH}
\Ind_K^G(\Ind_H^K(\Res_K^H( \Res_G^K(x))))\\
&\\
& =  \sum_{K\in \cF} b_K\sum_{H \subseteq K} c_{KH}
\Ind_H^G( \Res_G^H(x))
\end{array}
$$
We now define $a_H = \sum_{K \in \cF} b_K c_{KH}$,
and the formula becomes
$$x = \sum_{H\in \cF} a_H \Ind_H^G( \Res_G^H(x))\ .$$
\end{proof}

\section{Product formulas in surgery}
The existence of free actions on products of spheres will
be established by evaluating the surgery obstruction of a product
of degree one normal maps. We need to recall the product formulas in algebraic surgery theory due to Ranicki \cite[\S 8]{ra10}, \cite[\S 8]{ra11}. 
\begin{theorem}[{\cite[Prop. 8.1]{ra10}}]
Let $(A, \epsilon)$ and $(B, \eta)$ be rings with involution.
There are natural products in the symmetric and quadratic 
$L$-groups 
$$\begin{array}{lll}
\otimes:&L^n(A,\epsilon)\otimes_{\bZ} L^m(B,\eta)&\to
L^{m+n}(A\otimes_{\bZ} B, \epsilon\otimes \eta)\\
&&\\
\otimes:&L^n(A,\epsilon)\otimes_{\bZ} L_m(B,\eta)&\to
L_{m+n}(A\otimes_{\bZ} B, \epsilon\otimes \eta)\\
&&\\
\otimes:&L_n(A,\epsilon)\otimes_{\bZ} L^m(B,\eta)&\to
L_{m+n}(A\otimes_{\bZ} B, \epsilon\otimes \eta)\\
&&\\
\otimes:&L_n(A,\epsilon)\otimes_{\bZ} L_m(B,\eta)&\to
L_{m+n}(A\otimes_{\bZ} B, \epsilon\otimes \eta)
\end{array}
$$
\end{theorem}
These product formulas appear in the formula for  the surgery
obstruction of the product of normal maps.
Let $(f, b)\colon M \to X$ be a degree one normal map, where $M$ is a closed smooth or topological manifold of dimension $n$, 
$X$ is a finite Poincar\'e complex of dimension $n$, and $b\colon \nu_M \to \xi$ is a vector bundle map covering $f$. 
The symmetric signature $\sigma^*(X)$ is a cobordism invariant
of $X$  in the symmetric $L$-group $L^n(\bZ\pi_1(X))$, and
the quadratic signature $\sigma_*(f)$ of the normal map is a 
normal cobordism invariant lying in the quadratic $L$-group
$L_n(\bZ\pi_1(X))$. 
If $n\geq 5$, then $\sigma_*(f,b) = 0$ if and only if the normal map $(f,b)$ is normally cobordant to a homotopy equivalence \cite[\S 7]{ra11}.
These invariants are defined in \cite{ra11} by enriching the chain complexes $C(\widetilde M)$, $C(\widetilde X)$ of the universal covering spaces with additional structures arising from Poincar\'e duality. We will suppress mention of the reference maps $c_X\colon X \to K(\pi_1(X),1)$ and the orientation
characters, as well as the $K$-theory torsion
decorations (see \cite[\S 9]{ra10} and 
\cite[\S 6]{hrt1}).
\begin{theorem}[{\cite[Prop. 8.1(ii)]{ra11}}]
\label{thm: productformula}
Let $(f,b)\colon M \to X$ and $(g,c)\colon N \to Y$ be degree one
normal maps, with $\dim X = n$ and $\dim Y = m$.
Then for the product $(f\times g, b\times c)\colon M\times N \to X\times Y$ degree one normal map, the surgery obstruction
$$\sigma_*(f\times g,b\times c) = 
\sigma_*(f,b)\otimes \sigma_*(g,c) +
\sigma^*(X)\otimes \sigma_*(g,c)  +
 \sigma_*(f,b)\otimes \sigma^*(Y)$$
 as an element in $ L_{m+n}(\bZ[\pi_1(X\times Y)])$.
\end{theorem}
The naturality of the Ranicki product formulas 
(see \cite[Prop. 6.3]{hrt1}) will be used to reduce the
computation of surgery obstructions to hyperelementary groups.
Since we will be dealing with subgroups and finite coverings, it will be
convenient to introduce the  following notation. Let $X(G)$ denote
a finite Poincar\'e complex with fundamental group $G$, and
let $X(H)$ denote its covering space with fundamental group 
$H \subseteq G$. We have a similar notation 
$(f(H), b(H))\colon M(H) \to X(H)$ for coverings of
degree one normal maps. There are two basic operations
$$\Res_G^H\colon L_n(\bZ G) \to L_n(\bZ H)$$
and 
$$\Ind_H^G\colon L_n(\bZ H) \to L_n(\bZ G)$$
associated to subgroups, where the first is defined only for subgroups
of finite index. These give the $L$-groups a natural Mackey functor structure (see \cite[\S 1B]{htw3}, \cite[\S 5]{hrt1}). The symmetric and quadratic signatures behave
nicely under induction and restriction. For restriction, we can just consider $\sigma^*(X(G)) = (C(X(e)), \varphi_X)$ as a symmetric structure over $\bZ H$ by restricting the $G$-action to $H$.
 Here $e\in G$ denotes the identity element, so $X(e) = \widetilde X$.

For induction from subgroups, the induced symmetric structure is defined on the $G$-covering $G\times_H X(e) \to X(H)$ associated to the composition $X(H) \to K(H,1) \to K(G,1)$ of reference maps. This gives a symmetric structure denoted $\sigma^*(G\times_H X(e))$
on the $\bZ G$-chain complex $C(G\times_H X(e))$ (compare \cite[p. 196]{lee1}). 
Similarly, we have an \emph{induced} quadratic signature denoted
$$\sigma_*(G\times_H(f(e),b(e))\colon G\times_H M(e) \to G\times_H X(e))$$
 for a degree one normal map associated to
a covering $X(H)\to X(G)$. 
  In summary:
\begin{lemma}\label{lem: coverings}
Let $(f(G),b(G))\colon M(G) \to X(G)$ be an $n$-dimensional degree one normal map, with $G$ a finite group. For any subgroup
$H\subseteq G$ the following formulas hold:
\begin{enumerate}
\item \ $\Res_G^H(\sigma^*(X(G))) = \sigma^*(X(H))$.
\smallskip
\item \ $\Ind_H^G(\sigma^*(X(H))) = \sigma^*(G\times_H X(e))$.
\smallskip
\item \ $\Res_G^H(\sigma_*(f(G),b(G)))= \sigma_*(f(H),b(H))$.
\smallskip
\item \  $\Ind_H^G(\sigma_*(f(H),b(H)))= \sigma_*(G\times_H(f(e),b(e)))$.
\end{enumerate}
\end{lemma}
We now combine these formulas with Theorem \ref{thm: computable}.
\begin{theorem}\label{thm: reduction}
Let $(f,b)\colon M(G_1) \to X(G_1)$ and $(g,c)\colon N(G_2) \to Y(G_1)$ be degree one normal maps, with finite fundamental groups $G_1$ and $G_2$ respectively.
Then the quadratic signature $\sigma_*(f\times g,b\times c)=0$  if the products:
\smallskip
\begin{enumerate}
\item $\sigma_*(f(H_1),b(H_1))\otimes \sigma_*(g(H_2),c(H_2))=0$
\smallskip
\item $\sigma_*(f(H_1),b(H_1))\otimes\sigma^*(Y(H_2))=0$
\smallskip
\item $\sigma^*(X(H_1))\otimes \sigma_*(g(H_2),c(H_2))=0$
\end{enumerate}
\smallskip
for  all $2$-hyperelementary subgroups $H_1\subseteq G_1$
and $H_2\subseteq G_2$.
\end{theorem}
\begin{proof}
The $L$-groups needed for the product formulas are all 
$2$-hyperelementary computable after $2$-localization
(we will assume this is done, but not add to the notation). Moreover, Wall proved that $2$-localization is an injection for $L$-groups of finite groups, so the geometric surgery obstruction is $2$-locally detected.
Notice that by Lemma \ref{lem: coverings}, for any inclusion $H \subset G$ of finite groups,
$$\Ind_H^G(\Res_G^H\sigma^*(X(G))) = \Ind_H^G\sigma^*(X(H)),$$
and similarly for normal maps. 
By Theorem \ref{thm: computable},
we may write the symmetric
signature of $X(G)$ as an integral linear combination
of the images under induction from $2$-hyperelementary
subgroups
of the symmetric signatures of the coverings $X(H)$, and similarly
for the quadratic signatures of normal maps.

By substituting these expressions for $ G_1$ and $ G_2$  into each individual term of the product formula given in Theorem \ref{thm: productformula}, and
 using the
naturality of the product pairings under induction, 
we may write the quadratic signature
 $\sigma_*(f\times g,b\times c)\in L_{n+m}(\bZ[G_1\times G_2])$ as an integral linear combination of the 
 products listed,  for all $2$-hyperelementary
subgroups $H_1\subseteq G_1$ and $H_2\subseteq G_2$.
\end{proof}
 
\section{Periodic groups and normal Invariants}
The foundational work of Swan \cite{swan1} provides an $n$-dimensional, finitely dominated,
Poincar\'e complex $X=X(G)$ with fundamental group $G$,  for every finite group $G$ with 
periodic cohomology, and every sufficiently large $n\geq 3$ such that $n+1$ is a multiple of the period of $G$. It  has the property that the universal
covering $\widetilde X =X(e) \simeq S^n$. We will call these Swan complexes. 
The topological spherical
space form problem is to decide which Swan complexes are
homotopy equivalent to a closed   topological  manifold.
 
In  \cite[Thm 2.2]{thomas-wall1} it was shown the equivalence classes of $(G,n)$-\emph{polarized} Swan complexes (i.e. fixing an identification
of the fundamental group with $G$ and a homotopy equivalence of the universal covering
with $S^n$) correspond bijectively with generators of $H^{n+1}(G,\bZ)$. There is a slight difference in our notation from
that of \cite{madsen-thomas-wall1} since our groups are acting on
$S^n$ and not $S^{n-1}$. Recall that there is a classification
of periodic groups into types I-VI, and that some (but not all) of the periodic groups admit fixed-point free orthogonal representations. 
The quotients of $S^n$ by a fixed-point free orthogonal action
are called \emph{orthogonal spherical space forms}.
 
Here are the basic technical results of 
\cite[{\S\S 2-3}]{madsen-thomas-wall1}
and \cite[Thm. 2]{madsen-thomas-wall2}.
\begin{theorem}[{Madsen-Thomas-Wall}]\label{thm: mtw1}
Let $G$ be a finite group with periodic cohomology.
There exist finite $(G,n)$-polarized Swan complexes
$X=X(G)$ such that the covering spaces $X(H)$ 
are homotopy equivalent to  closed
manifolds, 
for each $H\subseteq G$ which has 
a fixed-point free orthogonal representation.
In addition, $X$ has a
smooth normal invariant which restricts to the normal
invariant of an orthogonal spherical space form for
the $2$-Sylow covering $X(G_2)$.
\end{theorem}
\begin{proof}
The first statement is Lemma 2.1 of \cite{madsen-thomas-wall1}, and the existence of smooth normal invariants is established in Theorem 3.1 of \cite{madsen-thomas-wall1}. The existence of a smooth normal invariant which restricts to that of an orthogonal spherical space form over the $2$-Sylow covering was established in \cite[Theorem 3.10]{madsen1}. This result is also contained in the proof of Theorem 2 in \cite{madsen-thomas-wall2}: the argument is based on Lemmas 3.2 and 3.3 of \cite{madsen-thomas-wall1}, which do not assume the $2p$-conditions, and the concluding step is on
\cite[p. 380]{madsen-thomas-wall1}. 
\end{proof}
 
We now recall the classification of $2$-hyperelementary
periodic groups. These are of type I or II, and those of type
I are semi-direct products of the form 
$$1 \to \cy{m} \to G \to  \cy{2^k} \to 1,$$
where $m$ is odd, and   $t\colon \cy{2^k} \to Aut(\cy{m})$ 
is the twisting map.
The type II groups are semi-direct products
$$1 \to \cy{m} \to G \to Q({2^\ell}) \to 1$$
 of odd order cyclic groups
by quaternion $2$-groups $Q(2^\ell)$.
\begin{lemma}
A $2$-hyperelementary periodic group $G$  has a non-central
element of order two if and only if $G$ is type I and $\ker t =\{1\}$.
\end{lemma}
We will call these $2$-hyperelementary type I groups with $\ker t =\{1\}$ the
\emph{generalized dihedral groups} and denote them 
$G=D(m,2^k)$. The ordinary dihedral groups $D_m$ of order
$2m$ are listed as $D(m,2)$
in this notation. Then a periodic group $G$ satisfies
the $2p$-conditions (for all primes $p$) if and only if
$G$ does not contain any (generalized) dihedral subgroups.
Note that the generalized dihedral groups are the only $2$-hyperelementary periodic groups which do \emph{not} admit a fixed-point free orthogonal representation.
 
The main result of \cite{madsen-thomas-wall1}, that a finite
group $G$ can act freely on some sphere $S^n$ if and only if 
$G$ satisfies the $p^2$ and $2p$ conditions for all primes $p$,
follows immediately from the results above. Here is a re-formulation
of the final surgery step.
\begin{corollary}[{Madsen-Thomas-Wall}]\label{cor: mtw2}
Let  $G$ be a finite group with periodic cohomology.
There exists a finite $(G,n)$-polarized Swan complex $X$, 
with $n\equiv \newmod{3}{4}$, 
and a smooth degree one normal
map $(f,b)\colon M \to X$ such that the quadratic signature
$\sigma_*(f,b) = 0 \in L_n(\bZ G)$ provided that
$G$ contains no (generalized) dihedral subgroups.
\end{corollary}
 \begin{proof}
We include a summary of their proof.
 If $G$ has no generalized dihedral subgroups, then \emph{every} $2$-hyperelementary subgroup of $G$ admits a fixed-point free orthogonal representation. One then starts with a finite $(G,n)$-polarized Swan complex $X$ as in Theorem \ref{thm: mtw1}, and the smooth degree one normal map $(f,b)\colon M \to X$ given by the special choice of normal invariant.
 Since $2$-localization gives an injective map on the surgery obstruction groups, Dress induction implies that
 $\sigma_*(f,b) = 0 \in L_n(\bZ G)$ provided that
 $$\sigma_*(f(H),b(H)) = 0 \in L_n(\bZ H)$$ for all $2$-hyperelementary
 subgroups $H\subset G$. However, each $X(H)$ is homotopy equivalent to a closed smooth manifold, and
surgery obstructions of normal maps between closed manifolds are detected by restriction to the $2$-Sylow subgroup $G_2 \subset G$. 
Therefore  $\sigma_*(f,b) = 0$.
This argument is explained in detail in \cite[\S 4]{madsen-thomas-wall1}.
\end{proof}
\begin{remark} The main result of \cite{madsen-thomas-wall1} gives free smooth actions on homotopy spheres, but not necessarily actions on the \emph{standard} sphere. It doesn't use the full strength of the special choice of smooth normal invariant provided by Theorem \ref{thm: mtw1}.
In  \cite[Theorem A]{madsen1}, Madsen showed that every group $G$ satisfying the $p^2$ and $2p$ conditions, for all primes $p$, acts freely and smoothly on the standard sphere. The argument follows from the naturality of the surgery exact sequence under the covering
$X(e) \to X$, and the surjectivity of the restriction map $L_{n+1}(\bZ G)\to L_{n+1}(\bZ)$. The special choice of normal invariant ensures that the smooth structure on $X(G)$ is covered by the standard sphere.
\end{remark}
\section{Surgery obstruction groups}
We will need some information about the surgery obstruction groups
$L_3(\bZ G)$ and $L_2(\bZ G)$, for $G$ a generalized dihedral group.
Here we must be precise about the $K$-theory decorations: since
every Swan complex $X$ is \emph{weakly-simple}, we can evaluate
the surgery obstruction $\sigma_*(f,b)$ of a degree one normal map
$(f,b) \colon M \to X$ in $L_n^\prime(\bZ G)$, where $L^\prime$ denotes
the weakly-simple obstruction groups (see \cite[\S 1]{htaylor2})
with allowable torsions in the subgroup $SK_1(\bZ G)\subset
Wh(\bZ G)$. Our Swan complexes will all have dimensions
$n\equiv \newmod{3}{4}$, so $\sigma_*(f,b)\in L^\prime_3(\bZ G)$.
The obstruction for the product of  two such problems will  be evaluated in 
$L^\prime_2(\bZ[G_1\times G_2])$, since $n + m \equiv \newmod{2}{4}$.
\begin{lemma}\label{lem: dihedral}
 Let $G$ be a generalized dihedral group,
and let $(f,b)\colon M \to X$ be a degree one normal map to
a finite $(G,n)$-polarized Swan complex.
\begin{enumerate}
\item The image of $\sigma_*(f,b)$ in $L_3^h(\Zhat_2 G)$ is
non-zero, but its image in $L_3^p(\Zhat_2 G)$ is zero.
\item If $G_1$ and $G_2$ are generalized dihedral groups,
the natural map $$L^\prime_2(\bZ [G_1\times G_2])
 \to L_2^p(\Zhat_2 [G_1\times G_2])$$ is
an injection.
\end{enumerate}
\end{lemma}
\begin{proof}
For the first statement we refer to \cite{lee1} or
\cite{jdavis1} for the result that the image of $\sigma_*(f,b)$ in $L_3^h(\Zhat_2 G)$ is detected by semi-characteristics, and
that these vanish if and only if $G$ satisfies the $2p$-conditions.
On the other hand, $L_3^p(\Zhat_2 G)=0$ since by reduction modulo
the radical we obtain the odd $L^p$-groups of a semisimple ring,
and these all vanish \cite{ra6}.

For the second part we start with 
the computation of $L^\prime_2(\bZ G)$.
The case $G = D(m,2)$ is done in \cite[\S 13]{htaylor2} and we follow
the outline given there. Recall that there is a natural splitting
(see \cite[\S 5]{htaylor1}):
$$L^\prime_2(\bZ G) = \bigoplus_{d\mid m}L^\prime_2(\bZ G)(d)$$
induced by idempotents in the $2$-localized Burnside ring. If 
$d\mid m$, then the restriction map $\Res_G^H$ is an isomorphism
on the $d$-component, where $H=N_G(\cy{d})$ is the subgroup
of $G$ whose odd-order part is $\cy d$. The splitting is compatible
with change of coefficients in $L$-theory.
\begin{lemma} If $G$ is a generalized dihedral group,  
 then $L^\prime_*(\bZ G)(d)$ is torsion-free for all $d\mid m$ such that
$-1 \in \text{Aut}(\cy d)$ is not contained in $\Image t$.
\end{lemma}
\begin{proof}
We combine \cite[Thm. 8.3]{htaylor2} with a 
calculation  $H^*(\wh'(\Zhat_2 G)(d))=0$, similar to that of \cite[\S 10]{htaylor2}, to obtain $L^\prime_*(\Zhat_2 G)(d)=0$. The result
that  $L^\prime_*(\bZ G)(d)$ is torsion-free now follows from \cite[Table 14.15]{htaylor2}.
\end{proof}
We assume from now on that $-1\in \Image t$ for all generalized
dihedral groups under consideration, since
 our $2$-torsion surgery obstructions clearly vanish
under restriction to the $d$-components where $L^\prime_*(\bZ G)(d)$ is torsion-free.
Under this assumption, 
 the various summands of $\bQ G$ all have type $O$, 
so from \cite[Table 14.12]{htaylor2} the natural map
$L^\prime_2(\bZ G) \to L^\prime_2(\Zhat_2 G)$
is an injection. From \cite[\S 9]{htaylor2} we check that the kernel
of the map
$$\psi_2\colon L^\prime_2(\Zhat_2 G)
 \to L^\prime_2(\bZ G \to \Zhat_2 G)$$
injects into $L^h_2(\Zhat_2 G)$. The basic point is that  $L^\prime_2(\Zhat_2 G)$ is a direct sum of terms of the
form $H^1(\hat A_2^{\times})\oplus g_2\cdot \cy 2$, where 
the summand $g_2\cdot \cy 2$ injects into $L^p_2(\Zhat_2 G)$,
and the summand $H^1(\hat A_2^{\times})$ injects into
$H^1((\hat A_2\otimes \bQ)^{\times})$ under the map $\psi_2$.

We now have a comparison sequence
$$\dots \to H^0(\tilde K_0(\Zhat_2 G)) \to 
L^h_2(\Zhat_2 G) \to L^p_2(\Zhat_2 G)\dots $$
where (by convention) we use the involution $[P] \mapsto -[P*]$
on $\tilde K_0(\Zhat_2 G))$. However, 
 $\tilde K_0(\Zhat_2 G))$ is free abelian,
with induced involution $-1$, so $H^0(\tilde K_0(\Zhat_2 G))=0$.

The steps for computing $L^\prime_2(\bZ [G_1\times G_2])$,
if $G_1$ and $G_2$ are generalized dihedral, are the same.
The summands of $\bQ [G_1\times G_2]$ all have type $O$,
and we get an injection into $L^\prime_2(\Zhat_2 [G_1\times G_2])$.
The same arguments as above show that
$\ker \psi_2$
injects into $L^p_2(\Zhat_2 [G_1\times G_2])$. 
\end{proof}
We also need one computation in the non-oriented case.
Let $G\times {\cy 2}^{-}$ denote the group $G\times \cy 2$ with
the non-trivial orientation character $w\colon
G\times \cy 2 \to \{\pm 1\}$ obtained by projecting onto
the second factor.

\begin{lemma}\label{lem: nonoriented}
Let $G=D(m, 2^k)$ be a generalized dihedral group. Then the 
natural map
$L_1^\prime(\bZ [G\times {\cy 2}^{-}])
\to L_1^h(\Zhat_2[G\times {\cy 2}^{-}])$
is an isomorphism.
\end{lemma}
\begin{proof}
The existence of a central element $T$ of order two with
$w(T)=-1$ implies that there are no involution invariant
simple factors of the rational group algebra. In other words,
we have type $GL$ for all factors,  so the relative
groups $L^\prime_*(\bZ G \to \Zhat_2 G)=0$, and 
$L_1^\prime(\bZ[G\times {\cy 2}^{-}])
 \cong L_1^\prime(\Zhat_2[G\times {\cy 2}^{-}])$. 
By Morita theory 
(see \cite[\S 8]{hrt1} for references), there is  an isomorphism
$$L_1^\prime(\bZ[G\times {\cy 2}^{-}])\cong \bigoplus_{d\mid m}
L_1^\prime(A(d)[{\cy 2}^{-}])$$
where $A(d) := \Zhat_2[\zeta_d]^{\sigma}$
is an unramified $2$-ring (with trivial involution) since
$\zeta_d$ is a primitive odd order  root of unity. The subgroup
$\sigma:=t^{-1}\langle 2\rangle\subseteq \cy{2^k}$ stabilizes each dyadic prime
in $\Zhat_2[\zeta_d]$, so acts as Galois automorphisms on the
ring of integers.
Now by \cite[Prop. 3.2.1]{wall-VI}, the Tate cohomology groups
$H^*(Wh^\prime(A(d)[{\cy 2}^{-}]) = 0$ (note the we are in the
\emph{exceptional case} mentioned in the paragraph immediately
following Prop. 3.2.1). By the Rothenberg sequence, it follows that
$L_1^\prime(\Zhat_2[G\times {\cy 2}^{-}])
\cong L_1^h(\Zhat_2[G\times {\cy 2}^{-}])$.
\end{proof}
In the next section, we will need to compute an $S^1$-transfer
on $L$-groups in one particular situation. If $S^1\to E \to B$
is an $S^1$ fibration,  let $\phi\colon \pi_1(B) \to \{\pm 1\}$ represent its first Stiefel-Whitney class. For any orientation charater 
$w\colon \pi_1(B) \to \{\pm 1\}$, we have a transfer map
$\Omega_n(B,w) \to \Omega_{n+1}(E,w\phi)$
defined on the the singular bordism groups of closed manifolds whose
orientation character is pulled back from $w$ or $w\phi$ respectively. 
If
$\bZ \to G \to \pi\to 1$ is the fundamental group sequence of the bundle,
then Wall \cite[Thm. 11.6]{wallbook} defined a transfer map
$L_n(\bZ \pi,w) \to L_{n+2}( \bZ G \to \bZ\pi,w\phi)$ and embedded it
in a long exact sequence with relative groups $LS_n(\Phi) \equiv LN_n(G \to \pi)$.
The \emph{$S^1$-bundle transfer} associated to the fibration is
the composition
$$L_n(\bZ \pi,w) \to L_{n+2}( \bZ G \to \bZ\pi,w\phi) \to L_{n+1}(\bZ G,w\phi)\ . $$
The product formulas are natural with respect to $S^1$-bundle transfers.
\begin{lemma}\label{lem: nattrans}
If $S^1 \to E \to B$ is a fibration with fundamental group sequence
$\,\bZ \to G_2 \to \pi\to 1$, then for any group $G_1$ the following
diagram commutes
$$\xymatrix{L_n(\bZ G_1) \otimes \Omega_{m-1}(B,w)\ar[r]\ar[d]&
L_{m+n-1}(\bZ[G_1\times \pi], 1\times w)\ar[d]\cr
L_n(\bZ G_1) \otimes \Omega_{m}(E,w\phi)\ar[r]&
L_{m+n}(\bZ[G_1\times G_2], 1\times w\phi)
}$$
where the vertical maps are $S^1$-bundle transfers and the horizontal
maps are obtained by products with the symmetric signature.
\end{lemma}
\begin{proof}
Let $(f,b)\colon M \to X$ be degree one normal map
 with fundamental group
$G_1$, and suppose that $S^1\to E \to B$ is an $S^1$-bundle
with the given fundamental group  and orientation data.
Elements of the cobordism group
$\Omega_{m-1}(B,w)$ are represented by 
pairs $(Y, h)$, 
where $Y$ is an $(m-1)$-dimensional closed manifold and $h\colon Y \to B$ is a reference map, with $w_1(Y) = h^*(w)$.  In this case, 
$\sigma_*(f,b)\otimes \sigma^*(Y,h)
=h_*(\sigma_*(f\times id, b\times id))$. In other words, the product
is represented by the quadratic signature associated to the degree
one  map  $f\times id_Y\colon M\times Y \to X \times Y$. Then we can define the $S^1$-bundle transfer by
pulling back the bundle $p\colon E\to B$ to $p^*E \to Y$, and
taking the quadratic signature of the $(m+n)$-dimensional degree one map
$f\times id\colon M \times p^*E \to X \times p^*E$.
\end{proof}

In \cite[p. 808]{ra12} the relative groups (in the 
quadratic case) were given 
an algebraic description 
$$LN_n(G \to \pi) = L_n(\bZ G, \beta, u)$$
where $(\beta, u)$ is an anti-structure 
on the group ring $\bZ G$ (see \cite[\S 1]{htaylor2}). 
If the base $B$ is a non-orientable
Poincar\'e space, but the total space $E$ is orientable, let
 $w\colon \pi \to \{\pm 1\}$ denote the orientation character for 
the base, and $\langle t\rangle\subseteq G$
denote the subgroup of $\pi_1(E)$ generated by $\pi_1(S^1)$.
We define $w\colon G \to \{\pm 1\}$ by composition.
Under these assumptions, $\beta(g) = g^{-1}$ if $w(g) = +1$ and
$\beta(g) = - g^{-1}t$ if $w(g) = -1$ (see \cite[p. 805]{ra12}).
\begin{example} Let $Q=Q(2^\ell)$ be a generalized quaternion
group. There
 exists a non-orientable fibre bundle
$$S^1 \to S^3/Q(2^\ell) \to RP^2$$
with 
 fundamental group sequence
$\bZ \to Q(2^\ell) \to \cy 2\to 1$, and orientable total space. The image of $\pi_1(S^1)$
 is a cyclic subgroup $\cy {2^{\ell -1}}\subseteq Q(2^\ell)$.
If $X =X(G)$ is a Swan complex, then the fibration
$$S^1 \to X \times S^3/Q(2^\ell)
\to X\times RP^2$$ realizes the fundamental group
sequence $\,\bZ \to G \times Q \to G \times \cy 2\to 1$.
\end{example}
\begin{corollary}\label{cor: transfer}
 Let $G$ be a generalized dihedral group, 
and $Q = Q(2^\ell)$ a generalized quaternion group.
The $S^1$-bundle transfer
$$L_1^\prime(\bZ [G\times {\cy 2}^{-}]) \to 
L_2^\prime(\bZ [G\times Q])$$
is zero.
\end{corollary}
\begin{proof}
The first map $L_1^\prime(\bZ [G\times {\cy 2}^{-}])
\to L_3^\prime(\bZ [G\times Q \to G\times {\cy 2}])$ in the definition
of the $S^1$-transfer is already zero. We will check this by  computing the previous term in the long exact sequence
$$\dots \to LN_1(G\times Q \to G\times \cy 2) \to L_1^\prime(\bZ [G\times {\cy 2}^{-}])
\to L_3^\prime(\bZ [G\times Q \to G\times {\cy 2}])$$
In this case, $LN_1(G\times Q \to G\times \cy 2) = 
L_1^\prime(\bZ[G\times Q, \beta, u])$ and we are computing
the map induced by projection
$$L_1^\prime(\bZ[G\times Q, \beta, u]) \to
 L_1^\prime(\bZ[G\times {\cy 2}^{-} ])$$
But we know from Lemma \ref{lem: nonoriented}
that $L_1^\prime(\bZ[G\times {\cy 2}^{-} ]) \cong 
L_1^h(\Zhat_2[G\times {\cy 2}^{-} ])$, so it is enough to
show that projection induces an isomorphism on the $2$-adic
$L$-groups. But this is clear since $Q(2^\ell)$ is a finite
$2$-group, and reduction modulo the radical induces an
isomorphism. Therefore
$L_1^h(\Zhat_2[G\times Q, \beta, u])\cong 
L_1^h(\Zhat_2[G\times {\cy 2}^{-} ])$.
\end{proof}

\section{The proof of Theorem A}
Let $G_1$ and $G_2$ be finite groups with periodic cohomology.
We choose finite polarized Swan complexes 
$X=X(G_1)$ and $Y = Y(G_2)$, as given in Theorem \ref{thm: mtw1},
of dimensions $n\equiv \newmod{3}{4}$ and $m\equiv\newmod{3}{4}$.
In particular, the covering spaces $X(H_1)$ and $Y(H_2)$
are homotopy equivalent to  closed
manifolds, 
for each $H_1\subseteq G_1$ and each $H_2\subseteq G_2$
which has
a fixed-point free orthogonal representation.
By Theorem \ref{cor: mtw2}, 
these Poincar\'e complexes also admit smooth degree one normal
maps $(f,b)\colon M \to X$ and $(g,c)\colon N \to Y$
such that the quadratic signatures 
$\sigma_*(f,b) = 0 \in L_n(\bZ G_1)$ (respectively
$\sigma_*(g,c) = 0 \in L_m(\bZ G_2)$ if and only if
$G_1$ (respectively $G_2$) contains no generalized dihedral subgroups. We remark that, by taking joins, it is possible to arrange
for $n=m$, as required in the statement of Theorem A. Our construction
actually gives free actions on $S^n\times S^m$ for any
$n$ and $m$ arising from the Swan complexes above.

\smallskip
We will now compute the product formula
(Theorem \ref{thm: productformula}) to show that the quadratic
signature
$\sigma_*(f\times g, b\times c) = 0 \in L_{m+n}(\bZ[G_1\times G_2])$.
It will then follow by surgery theory, that
there is a smooth, closed manifold $W \simeq X\times Y$ such
that the universal covering $\widetilde W$ is homotopy equivalent to $S^n\times S^m$. Since $n+m \equiv \newmod{2}{4}$,
the surgery exact sequence \cite[\S 10]{wallbook} implies
that $\widetilde W$ is diffeomorphic to $S^n\times S^m$,
and we obtain a free action on a product of 
standard spheres.

Now we proceed with the calculation of the quadratic signature.
According to Theorem \ref{thm: reduction},
$\sigma_*(f\times g, b\times c) = 0$  provided that:
\smallskip
\begin{enumerate}
\item $\sigma_*(f(H_1),b(H_1))\otimes \sigma_*(g(H_2),c(H_2))=0$
\smallskip
\item $\sigma_*(f(H_1),b(H_1))\otimes\sigma^*(Y(H_2))=0$
\smallskip
\item $\sigma^*(X(H_1))\otimes \sigma_*(g(H_2),c(H_2))=0$
\end{enumerate}
\smallskip
for  all $2$-hyperelementary subgroups $H_1\subseteq G_1$
and $H_2\subseteq G_2$.
The calculation is now reduced to evaluating these products
in the  following three situations:
\smallskip
\begin{enumerate}
\item[A1.] $L^\prime_3(\bZ H_1)\otimes L_3^\prime (\bZ H_2) \to
 L_2^\prime(\bZ[H_1\times H_2])$,
where $H_1$ and $H_2$ are both generalized dihedral,
\smallskip
\item[A2.] $L^\prime_3(\bZ H_1)\otimes L^3 (\bZ H_2) \to
 L_2^\prime(\bZ[H_1\times H_2])$,
where $H_1$ and $H_2$ are both generalized dihedral,
\smallskip
\item[A3.] $L_3^\prime(\bZ H_1)\otimes L^3(\bZ H_2)
\to L_2^\prime(\bZ[H_1\times H_2])$,
where $H_1$ is generalized dihedral, and $H_2$ has
no dihedral subgroups.
\end{enumerate}
In addition, we may reverse the roles of $H_1$ and $H_2$ in A2 and A3 to cover all the possible cases. No new calculations are involved.
\subsection*{Cases A1 and A2}
By Lemma \ref{lem: dihedral}, 
 the image of $\sigma_*(f,b)$ in $L_3^h(\Zhat_2 H_1)$ is
non-zero, but its image in $L_3^p(\Zhat_2 H_1)$ is zero.
On the other hand, 
the natural map $L^\prime_2(\bZ [H_1\times H_2])
 \to L_2^p(\Zhat_2 [H_1\times H_2])$ is
an injection.
The vanishing of $\sigma_*(f,b)\otimes \sigma_*(g,c)$
in Case A1
follows from applying these facts to the diagram
$$
\xymatrix{L_3^h(\Zhat_2 H_1)\otimes L_3^h(\Zhat_2 H_2)
\ar[r]\ar[d]& L_2^h(\Zhat_2[H_1\times H_2])\ar[d]\cr
L_3^p(\Zhat_2 H_1)\otimes L_3^p(\Zhat_2 H_2)
\ar[r]& L_2^p(\Zhat_2[H_1\times H_2])
}
$$
For Case A2 we use a similar diagram with $L^3(\Zhat_2 H_2)$
instead of $L_3^h(\Zhat_2 H_2)$. We get  $\sigma_*(f,b)\otimes \sigma^*(Y)=0$, and similarly  $\sigma^*(X)\otimes\sigma_*(g,c)=0$
by reversing the roles of $H_1$ and $H_2$.

\subsection*{Case A3}
In this case, $H_2$ has no dihedral subgroups, so $X(H_2)$ is
homotopy equivalent to an orthogonal spherical space form.
Therefore $\sigma^*(Y)$ is the symmetric signature of a closed
manifold, and hence lies in the image 
$$\Omega_m(BH_2) \to L^m(\bZ H_2)$$
This subgroup has better induction properties: it is computable using the family of $2$-Sylow subgroups (recall that we have localized at
2). It follows that
$\sigma^*(Y)$ as an odd integral multiple of 
$\Ind_{K_2}^{H_2}(\sigma^*(Y(K_2)))$, where $K_2\subseteq H_2$ is
a $2$-Sylow subgroup. By naturality of the product formula, we may  assume from the beginning
that $H_2= \cy{2^\ell}$ or $H_2=Q(2^{\ell})$.

Since $\Omega_m(BH_2)\otimes \bZ_{(2)} \cong H_m(BH_2;\bZ_{(2)})\cong \cy{|H_2|}$, any orthogonal
space form $Y(H_2)$ is bordant to an (odd) multiple of a generator
of this homology group. Suppose that $m=4k+3$. Then we can represent a generator of $\Omega_m(BH_2)\otimes \bZ_{(2)}$
by a standard quotient $Y_0(H_2):=S^m/H_2$, constructed as the join
of $k+1$ copies of a free orthogonal action on $S^3$. By linearity of
the product formula, it is enough
to evaluate the product formula for $Y_0(H_2)$.

If $H_2=\cy{2^\ell}$,  then $Y_0=S^m/\cy{2^\ell}$ is the total space of an oriented
$S^1$-bundle:
$$ S^1 \to S^m/\cy{2^\ell} \to \bC P^{2k+1}$$
and $\bC P^{2k+1}$ is an oriented boundary. Therefore, the
product  $\sigma_*(f,b)\otimes \sigma^*(Y_0)$ is the 
$S^1$-bundle transfer of the product $\sigma_*(f,b)\otimes
\sigma^*( \bC P^{2k+1}) =0$.
\begin{remark}
An alternate argument could be given by computing the group
$L_2^\prime(\bZ[H_1\times H_2])$ using \cite[Prop.  3.2.3]{wall-VI}
as above.
\end{remark}

If $H_2=Q(2^{\ell})$, then $Y_0=S^m/Q(2^\ell)$ is the total space of a
non-orientable $S^1$-bundle
$$ S^1 \to S^m/Q(2^\ell) \to \bC P^{2k+1}/\cy{2}$$
where the base is the orbit space of $\bC P^{2k+1}$ under
a free, orientation-reversing $\cy 2$-action. In homogeneous
co-ordinates the action of the generator is given by
$$[z_0, z_1, \dots, z_{2k}, z_{2k+1}]\mapsto
[-\bar z_1, \bar z_0, \dots, -\bar z_{2k+1}, \bar z_{2k}]$$
where $\bar z$ denotes complex conjugation. There is a commutative
diagram of quotients
$$\xymatrix{ S^1\ar[r]\ar[d]&  S^m\ar[r]\ar[d] & \bC P^{2k+1}\ar[d]\cr
 S^1/\cy{2^{\ell-1}}\ar[r]& S^m/Q(2^\ell) \ar[r]& \bC P^{2k+1}/\cy{2}
}$$
arising from embedding the extension 
$1\to \cy{2^{\ell-1}} \to Q(2^\ell) \to \cy 2\to 1$ in the unit quaternions $S^3$, and taking
joins to get $S^{4k+3} = S^3 \ast \dots \ast S^3$.
Let $Z_0=\bC P^{2k+1}/\cy{2}$ denote this orbit space.

The product $\sigma_*(f,b)\otimes \sigma^*(Z_0)$ lies in
$L_1^\prime(\bZ [H_1\times {\cy 2}^{-}])$, where $H_1$
is a generalized dihedral group. 
By Corollary \ref{cor: transfer}, the 
 $S^1$-bundle transfer
$$L_1^\prime(\bZ [H_1\times {\cy 2}^{-}]) \to 
L_2^\prime(\bZ [H_1\times H_2])$$
is zero, so the product
$\sigma_*(f,b)\otimes \sigma^*(Y_0)=0$ by 
Lemma \ref{lem: nattrans}. This completes the proof of Theorem A.

\medskip
The
multiplicativity $\sigma^*(X \times Y) = \sigma^*(X) \otimes
\sigma^*(Y)$ of the symmetric signature
\cite[Prop. 8.1(i)]{ra11}, and the 
vanishing of pair-wise products established above, gives many more examples of the same kind:
\begin{theorem}\label{thm: moreex}
For $k >1$, $k\neq 3$, any product $G_1\times \dots \times G_k$ of finite
groups with periodic cohomology acts freely
and smoothly on $(S^n)^k$, for some $n\equiv \newmod{3}{4}$. 
\end{theorem}
\begin{remark}
If $k=3$, the same argument only shows that a product $G_1 \times G_2 \times G_3$ of periodic groups acts freely and smoothly either on $(S^n)^3$, or  possibly on
$(S^n)^3\, \sharp \, \Sigma ^{3n}$, where $\Sigma^{3n}$ denotes the Kervaire homotopy sphere in  dimension $3n\equiv \newmod{1}{4}$. If the Kervaire sphere is non-standard (e.g. when $3n \neq 2^{t} -3$, $t \geq 2$), then $(S^n)^3\, \sharp \, \Sigma ^{3n}$
is \emph{not} diffeomorphic to $(S^n)^3$. The smallest example is
$n=3$, and $\Sigma^9$ is non-standard. 
\end{remark}

\bibliographystyle{ih}
\providecommand{\bysame}{\leavevmode\hbox to3em{\hrulefill}\thinspace}
\providecommand{\MR}{\relax\ifhmode\unskip\space\fi MR }
\providecommand{\MRhref}[2]{%
  \href{http://www.ams.org/mathscinet-getitem?mr=#1}{#2}
}
\providecommand{\href}[2]{#2}

\end{document}